\documentclass[10pt]{amsart}
\usepackage{amsmath,amssymb,array,amsthm,epsfig}
\usepackage[T1]{fontenc}
\usepackage{amsfonts}
\usepackage{amsthm}
\usepackage{setspace}
\usepackage{tikz}
\usepackage[all,cmtip]{xy}
\usepackage{booktabs}
\usepackage{hyperref}
\usepackage{longtable}
\usepackage{geometry}
\usetikzlibrary{trees}

\usepackage[english]{babel}
\usepackage[utf8]{inputenc}
\usepackage{color}
\usepackage{mathrsfs}
\usepackage{latexsym}

 \makeindex

\usepackage{geometry}
\usepackage{marginnote}
\usepackage{appendix}

\usepackage{calligra}

\usepackage{color}
\usepackage[normalem]{ulem}

\usepackage{array}
\usepackage{makecell}
\usepackage{tabularx}
\usepackage[framemethod=tikz]{mdframed}

\newcommand{\nwc}{\newcommand}

\newcommand{\cal}{\mathcal}

\newcommand{\iz}{\mathcal{I}_{Z}}

\newcommand{\PP}{\mathbb{P}}

\nwc{\aaa}{\sF }
\nwc{\aab}{\bar{\mathfrak{a}}}
\nwc{\aal}{\sF '}
\nwc{\aap}{\sF _{P}}
\newcommand{\OO}{{\mathcal{O}}}

\nwc{\bbb}{\mathfrak{b}}
\nwc{\bbp}{\mathfrak{b}_{P}}

\nwc{\C}{\mathbb{C}}
\nwc{\cb}{\overline{C}}
\nwc{\ccc}{\mathcal{C}}
\nwc{\ch}{\widehat{C}}
\nwc{\cin}{\textbf{(v)}}
\nwc{\cl}{C'}
\nwc{\cp}{\mathcal{C}_{P}}
\nwc{\cpll}{\mathcal{C}_{P'}}
\nwc{\ct}{\tilde{C}}

\nwc{\dd}{\mathcal{D}}
\nwc{\ddd}{\mathfrak{d}}
\nwc{\ddl}{\mathcal{L}'}
\nwc{\dlp}{\delta_{P}}
\nwc{\doi}{\textbf{(ii)}}

\nwc{\ff}{\mathscr{F}}

\nwc{\G}{{\cal G}}
\nwc{\gon}{{\rm gon}}
\nwc{\gtl}{\tilde{g}}
\nwc{\gud}{g^{1}_{2}}
\nwc{\gtu}{g^{1}_{3}}

\nwc{\hhza}{H^{0}(C,\mathfrak{a})}
\nwc{\hua}{h^{1}(C,\mathfrak{a})}
\nwc{\hza}{h^{0}(C,\mathfrak{a})}

\def \d { \partial}

\nwc{\kk}{{\rm K}}

\nwc{\lbd}{\lambda}
\nwc{\lif}{L_{\infty}}

\nwc{\mm}{\mathfrak{m}}
\nwc{\mmp}{\mathfrak{m}_{P}}
\nwc{\mpd}{{\mathfrak{m}_{P}}^{2}}

\nwc{\N}{I\!\!N}
\nwc{\nn}{\mathbb{N}}

\nwc{\obp}{\overline{\mathcal{O}}_P}
\nwc{\ocbux}{\oo _{\bar{C}}\langle 1,x\rangle}
\nwc{\oclux}{\oo _{C'}\langle 1,x\rangle}
\nwc{\ocux}{\oo _{C}\langle 1,x\rangle}
\nwc{\ol}{\mathcal{O}'}
\nwc{\oma}{\Omega (\mathfrak{a})}
\nwc{\omo}{\Omega (\mathcal{O})}
\nwc{\oo}{\mathcal{O}}
\nwc{\ooh}{\widehat{\mathcal{O}}}
\nwc{\opc}{\mathcal{O}_{P,C}}
\nwc{\oph}{\widehat{\mathcal{O}}_{P}}
\nwc{\opl}{\mathcal{O}_{P}'}
\nwc{\oplc}{\mathcal{O}_{P,C}'}
\nwc{\opll}{\mathcal{O}_{P'}}
\nwc{\opt}{\tilde{\mathcal{O}}_{P}}
\nwc{\optt}{{\mathcal{O}}_{\tilde{P}}}
\nwc{\oq}{\mathcal{O}_{Q}}
\nwc{\oqt}{\tilde{\mathcal{O}}_{Q}}
\nwc{\ot}{\tilde{\mathcal{O}}}
\nwc{\overop}{\bar{\oo}_{P}}

\nwc{\pb}{\overline{P}}
\nwc{\pgmd}{\mathbb{P}^{g+2}}
\nwc{\pgmu}{\mathbb{P}^{g+1}}
\nwc{\pp}{\mathbb{P}}
\nwc{\prv}{\noindent\textbf{Proof}:}
\nwc{\pt}{\tilde{P}}
\nwc{\ptl}{\tilde{P}}
\nwc{\pum}{\mathbb{P}^{1}}
\nwc{\carta}{\mathfrak{U}}

\def\F{{\mathscr{F}}}

\def\sL{{\mathscr{L}}}
\renewcommand{\H}{\text{\rm H}}

\nwc{\Q}{\;\mbox{{\sf I}}\!\!\!Q}
\nwc{\qb}{\overline{Q}}
\nwc{\qtl}{\tilde{Q}}
\nwc{\qua}{\textbf{(iv)}}

\nwc{\R}{I\!\!R}

\nwc{\sep}{\beq\ast\ \ast\ \ast\enq}
\nwc{\spl}{{S_{P}}'}
\nwc{\spll}{S_{P'}}
\nwc{\ssp}{{\rm S}_{P}}
\nwc{\sss}{{\rm S}}
\nwc{\sys}{\mathcal{L}}

\nwc{\tre}{\textbf{(iii)}}

\nwc{\um}{\textbf{(i)}}

\nwc{\vlp}{\mathcal{V}_{\lambda,P}}
\nwc{\vpt}{v_{\ptl}}
\nwc{\vv}{\mathcal{W}}
\nwc{\vvp}{\mathcal{W}_{P}}
\nwc{\vpb}{v_{\overline{P}}}
\nwc{\vtxp}{\widetilde{V}_{x,P}}
\nwc{\vxp}{V_{x,P}}
\nwc{\vzp}{V_{z,P}}

\nwc{\wol}{\ww\cdot\mathcal{O}'}
\nwc{\wpn}{{\omega _{P}}^{n}}
\nwc{\wwt}{\widetilde{\omega}}
\nwc{\wwtp}{\widetilde{\omega}_P}
\nwc{\ww}{\omega}
\nwc{\wwp}{\omega _{P}}

\nwc{\zz}{\mathbb{Z}}

\def\sF{{\mathscr{F}}}

\input{xypic.tex} \xyoption{all}

\usepackage{float}
\setcellgapes{3pt}

\newtheorem{coro}{Corollary}[section]
\newtheorem{defi}[coro]{Definition}
\newtheorem{exem}[coro]{Example}

\newtheorem{lema}[coro]{Lemma}
\newtheorem{prop}[coro]{Proposition}

\newtheorem{teo}[coro]{Theorem}

\usepackage[utf8]{inputenc}

\usepackage{hyperref}
\hypersetup{
	colorlinks=true,
	linkcolor=blue,
	filecolor=magenta, 
	urlcolor=cyan,
}
\makeatletter
\@namedef{subjclassname@2020}{%
 \textup{2020} Mathematics Subject Classification}
\makeatother

\author{Maur\'icio Corr\^ea}
\address{\newline
\indent
Universit\`a degli Studi di Bari, 
Via E. Orabona 4, I-70125, Bari, Italy
\newline
\newline
\indent
UFMG\\
Avenida Ant\^onio Carlos, 6627\\
30161-970 Belo Horizonte\\ Brazil
}
\email{correajr@ufmg.br}



\title{Analytic varieties invariant  by   foliations and Pfaff systems }

\begin{document}

\begin{abstract}
 In this work 
we shall present a survey on problems and results on singular holomorphic
foliations and Pfaff systems on complex manifolds assuming that these objects possess invariant analytic
varieties. We will focus on recent results which have been motivated by classical works of Darboux,
Poincar\'e and Painlev\'e on the problem of algebraic integration of singular
polynomial differential equations. We present results on Poincar\'e and Painlev\'e problem of bounding the
degree and the genus of analytic varieties invariant by holomorphic foliations and Pfaff systems. We shall
discuss the general ideas of the theory of integrability characterizing the existence of meromorphic first
integrals for complex analytic Pfaff equations.  
\end{abstract}

\maketitle
\tableofcontents

\section{Introduction}
 
The study  of  regular distributions and foliations and their integral manifolds  was  motivated by the  classical work due to Pfaff \index{Pfaff}  \cite[Chapter III]{For}. The  study of singular  polynomial differential equations on
the complex plane  was investigated by Poincar\'e, Darboux, Painlev\'e.  
 Nowadays, this corresponds to the study of  singular  holomorphic foliations on complex projective
spaces.  G. Darboux presented  in \cite{Darboux} \index{Darboux} a theory on the existence of
first integrals for polynomial differential equations based on the
existence of sufficiently many invariant algebraic curves.
In \cite{Poincare} H. Poincar\'e \index{Poincar\'e}   observed that  it would be sufficient to
bound  the degree of algebraic solutions.  P.
Painlev\'e\index{Painlev\'e}, in  \cite{Painleve}, stated an integrability problem as follows:
 Is it possible to recognize the genus of the general
solution of an algebraic differential equation  \index{Algebraic differential equation} in two variables
which has a rational first integral? \  It is worth mentioning that  Painlev\'e  won  the  ``Grand Prix des Sciences Math\'ematiques'' of the French Academy of Sciences  by his important contributions in this subject, see \cite{Schlomiuk}.

These problems are  known as
 Poincar\'e's    and  Painlev\'e's   Problems. In \cite{Alcides-exemplo}
A. Lins Neto  constructed families of foliations on
$\mathbb{P}^2$, with fixed degree and local analytic
type of the singularities,  where foliations with rational first
integrals of arbitrarily large degree appear. In other words, such
families show that the  questions of Poincar\'e and Painlev\'e  have a
negative answer in general. However, one can obtain an affirmative
answer provided that some additional hypotheses are assumed.  

Several works,  such as  D. Cerveau and A. Lins Neto \cite{Cerv1},  M.
Carnicer \cite{Carnicer}, M. Soares \cite{Soares-Annals,Soares-Inventiones} and  Walcher \cite{Walcher}
 have stimulated the current interest in the investigation of the  Poincar\'e problem.
Many authors have been working on these problems and on some of its generalizations, see for instance
the papers by  M.  Soares \cite{Soares-Inventiones},   M. Brunella
 and  L.G.  Mendes \cite{BruMe}, E. Esteves  and  S. Kleiman \cite{Esteves-Kleiman-caracteristica}, V. 
Cavalier  and  D.  Lehmann \cite{Cavalier-Lehmann}. This problem  also has  been  considered in other  varieties such as: del Pezzo surfaces and K3 surfaces \cite{Zamora}, weighted projective spaces \cite{BCR, Correa-Soares-pesos}, multiprojective spaces  \cite{multiprojective},  toric varieties \cite{Miguel},  projective manifolds with Picard number equal to one  \cite{BruMe}, surfaces with trivial Picard group \cite{Ballico-Poincare},  varieties  over an
algebraically closed field of arbitrary characteristic \cite{Esteves-Kleiman-caracteristica}. It follows from a celebrated theorem by Jouanolou \cite{Jouanolou2} \index{Jouanolou} that foliations on projective plane, of degree at least 2,  with some invariant  algebraic curve are rare. We refer  the reader to  \cite{Coutinho-Pereira}  and references therein, where  Coutinho and Pereira provide  a generalization of Jouanolou's  result  for one dimensional
foliations and Pfaff fields  on  projective manifolds.

The problem of bounding the genus of an invariant curve in terms of
the degree of a one-dimensional foliation on
$\mathbb{P}^n$ was considered  for instance by
Campillo, Carnicer and de
 la Fuente \cite{Campillo-Carnicer-Fuente}. They showed  that, if  $C$ is a reduced curve which is invariant by a
 one-dimensional foliation  $\F$, of degree $d$, on  $\mathbb{P}^n$,  then
\begin{equation}\label{ccf}
\frac{2p_a(C)-2}{\deg(C)}\leq d-1+a,
\end{equation}
where $p_a(C)$ is the arithmetic genus \index{Arithmetic genus} of $C$ and $a$ is an integer
obtained from the concrete problem of imposing singularities
 to projective hypersurfaces. For instance, if $C$ has only nodal singularities \index{Nodal singularities} then $a=0$, and thus formula $(\ref{ccf})$
 follows from \cite{Fu}.
 Esteves and Kleiman in \cite{Esteves-Kleiman-caracteristica} have provided   bounds for the arithmetic genus of a  curve invariant by a foliation by curves, 
 which improve and extend some    results
of Campillo, Carnicer, and de la Fuente, and of du Plessis and Wall \cite{Plessis-Wall}.
  In  \cite{Correa-Jardim-genera}  we establish  upper bounds for the sectional genus of
Gorenstein varieties which are invariant by  Pfaff systems on projective spaces and in \cite{Ballico-Pfaff} Ballico gave   an extension of this result.

 The work of J.P.
Jouanolou in \cite{Jouanolou2}   gives  an  improvement and generalization
of the Darboux theory of integrability characterizing the existence
of rational first integrals for Pfaff equations  of codimension one on
$\mathbb{P}^n_{k}$, where $k$ is an algebraically closed field of
characteristic zero. Namely,  consider  $\omega\in \mathrm{H}^0(\mathbb{P}^n_{k},
\Omega_{\mathbb{P}^n_{k}}^1(d+1))$   a twisted $1$-form which induces a codimension one  Pfaff system.  Then by \cite[Theorem 3.3]{Jouanolou2}, we have that 
$\omega$  admits a rational  first integral if and only if it 
has an  infinite number of   invariant irreducible  hypersurfaces. More generally, 
Jouanolou proved in \cite{Jouanolou1} that on a complex compact manifold $X$
satisfying certain conditions on its Hodge-to-de Rham spectral
sequence, a Pfaff system $\omega\in \mathrm{H}^0(X,
\Omega_{X}^1\otimes\sL)$, where $\sL$ is a line
bundle, admits a meromorphic first integral if and only if  it admits 
an infinite number of invariant irreducible divisors. Moreover, if
$\omega$ does not admit a meromorphic first integral, then the
number of invariant irreducible divisors is at most
\begin{center}
$
\mathrm{dim}_{\mathbb{C}}(\mathrm{H}^0(X,\Omega_{X}^2\otimes \sL)/\omega\wedge\mathrm{H}^0(X,\Omega_{X}^1)) +\rho(X)+1,$
\end{center}
where $\rho(X)$ is the Picard number of $X$. Deschamps  provides a proof of this result for projective surfaces in his monograph  \cite{Deschamps} on  Bogomolov's work \cite{Bogomolov2} on  the boundedness of families of curves of fixed geometric genus on a   surface of general type with  positive Segre class.

E. Ghys \index{Ghys} in \cite{Ghys} drops  the cohomological  hypotheses given by Jouanolou showing
that this result is valid for all compact complex manifolds. 
More precisely, the result of
\cite{Ghys}   states that if a codimension one  Pfaff system 
$\omega \in H^0(X, \Omega^1_X \otimes \mathscr{L})$,
does not admit a meromorphic first integral, then the number of
invariant irreducible divisors must be smaller than
$$
  \mathrm{dim}_{\mathbb{C}}\left[\mathrm{H}^0(X,\Omega_{X}^2\otimes  \mathscr{L})/\omega\wedge\mathrm{H}^0(X,\Omega_{cl}^1)\right] + \mathrm{dim}_{\mathbb{C}} \mathrm{H}^1(X,\Omega_{cl}^1) + 2,
$$
where $\Omega_{cl}^1$ denotes the sheaf of closed holomorphic
$1$-forms on $X$.  In the case for holomorphic foliations of dimension one,   we presented in \cite{Correa-dim1}  a similar result which reads as follows:
let $X$ be as above and let $\F$ be a one-dimensional holomorphic foliation on $X$. If $\F$ admits at least
$$
\dim_{\mathbb{C}}\left[\mathrm{H}^0(X,T\F^*)/i_{v_{\F}}(\mathrm{H}^0(X,\Omega_{cl}^1))\right]+\dim_{\mathbb{C}} \mathrm{H}^1(X,\Omega^1_{cl})+\dim(X)
$$
invariant irreducible hypersurfaces, then $\F$ admits a meromorphic first integral. Here, $T\F$ is the  tangent  bundle of $\F$ and $i_{v_{\F}}(\cdot)$
denotes  the contraction by the local  vector fields   inducing the foliation $\F$. A higher codimensional version   of these results  has been proved in \cite{CMM}.  More precisely, we prove that 
if  $\F$ is a  Pfaff system, of  codimension $k$,  on a   compact  complex
manifold $X$,  defined by  $\omega\in H^{0}(X,
\Omega^{k}_X \otimes \mathscr{L})$, such that  $\F$ admits
\begin{equation}\label{G0}
\dim_{\mathbb{C}}\left[H^{0}(X,
\Omega^{k+1}_X \otimes \mathscr{L}) / \omega \wedge H^{0}(X,
\Omega_{cl}^1)\right] + \dim_{\mathbb{C}}H^{1}(X, \Omega_{cl}^{1}) + k + 1
\end{equation}
invariant irreducible analytic hypersurfaces, then $\F$ admits a meromorphic first integral. Versions of Darboux-Jouanolou theorem for polynomial differential forms have been provided in \cite{CMM2,Llibre-Zhang,Correa-campos-invariantes, Edileno-Sergio}. 

A discrete dynamical version of Darboux-Jouanolou's theorem was
 given   by S. Cantat.  He proved  in  \cite{Cantat} that if
there exist $N$ irreducible hypersurfaces  invariant  by an automorphism $f \colon X\to X$ with
$$
N\geq \dim(X)+ \mathrm{dim}_{\mathbb{C}}(\mathrm{H}^1(X,\Omega_{X}^1))
$$
then $f$ preserves a nontrivial meromorphic fibration.

M. Brunella and M. Nicolau in \cite{Brunella-Nicolau}  prove a version of the Darboux-Jouanolou theorem  for
 codimension one foliations  in positive characteristic. The authors have improved a previous result due to Kim in \cite{Kim}.  Also, in   \cite{Brunella-Nicolau} the authors  provide a Darboux-Jouanolou type theorem for non-singular
codimension one transversely holomorphic foliations on compact
manifolds. In \cite{cama} L. C\^amara has proved a high dimensional version of  Brunella and Nicolau   result for  transversely  holomorphic foliations.  Sc\'ardua in \cite{Bruno2} shows, under certain conditions, a  
Darboux-Jouanolou type  theorem for codimension one germs of holomorphic foliations on $(\mathbb{C}^n, 0)$. For   similar results in dimension   2 we refer to \cite{Camacho-Bruno,Bruno1}.  
A version of Darboux-Jouanolou theorem for singular toric varities was proved in   \cite{Correa-toricas} by using   Cox's homogeneous coordinate ring  \cite{Cox}\index{Cox's ring}.

\subsection*{Acknowledgments}
The author  is supported by CNPq grant numbers 202374/2018-1, 302075/2015-
1, 400821/2016-8; he is thankful to  Israel Vainsencher, Leonardo C\^amara, Miguel Pe\~na, Alan Muniz, Marcos Jardim,   Bruno Sc\'ardua and Marcio Soares  for useful comments. 
  
\section{Pfaff systems, distributions and foliations}

\subsection{Distributions and foliations}
Throughout this survey we denote by $X$ a   complex  manifold  of dimension $n$.

A  regular holomorphic foliation $\F$ of codimension $q$ on a complex manifold $X$  is given by an open covering $\{U_{\alpha}\}_{\alpha}$  of $X$, 
holomorphic  submersions $ f_{\alpha}:U_\alpha \to \mathbb{C}^q$  and  biholomorphisms 
$$g_{\alpha \beta}:f_\beta(U_{\alpha}\cap U_{\beta})\subset \mathbb{C}^q \to f_\alpha(U_{\alpha}\cap U_{\beta})\subset \mathbb{C}^q$$ 
 such that $f_\alpha = g_{\alpha \beta}\circ f_\beta$. This gives us a holomorphic vector bundle $T\F:=\cup_{\alpha} \ker(df_{\alpha})\subset TX $ called tangent bundle of $\F$. It follows from Frobenius theorem that a  holomorphic subvector bundle $E\subset TX$ (called distribution) of rank $n-q$   is a tangent bundle of a foliation of codimension $q$ if and only if $[E,E]\subset E$.

In this survey we will be interested in holomorphic  foliations and distributions with singularities.

\begin{defi}\index{Distributions}
A  singular  holomorphic  distribution $\F$  of codimension  $p$ on $X$ is a nonzero coherent subsheaf $T\F\subsetneq TX$
of generic rank $n-p$ which is  saturated, i.e.,
such that $TX / T\F$ is torsion free. 
\end{defi}

We have an exact sequence of sheaves 
$$
 0  \longrightarrow T\F  \longrightarrow TX \longrightarrow  TX / T\F:=N\sF  \longrightarrow 0.
$$
The sheaves $T\sF$ and $N\sF$ are called the \emph{tangent} and the \emph{normal} sheaves of $\sF$, respectively.
Now, by taking the double dual of the  $p$-th wedge product  of the inclusion
$$
N\sF^{*} \longrightarrow \Omega_X^1
$$
we get a map 
$$(\wedge^{p}N\sF^*)^{**} \longrightarrow \Omega^{p}_X.$$ 
Since $N\sF$ and $N\sF^*$ are   torsion-free,   it follows from  \cite[Proposition 5.6.10]{Kob} and \cite[Proposition 5.6.12]{Kob} that $(\wedge^{p}N\sF^*)^{**}\simeq \det(N\sF^*)\simeq \det(N\sF)^*$. This  gives rise to a nonzero twisted  holomorphic  $p$-form  $\omega_{\F}\in H^0(X, \Omega^{p}_X\otimes \det(N\sF)^{**})\simeq
 H^0(X, \Omega^{p}_X\otimes \det(N\sF))$.   The  twisted    $p$-form $\omega_{\F}$ is called  by the  {\it  Pfaff system}  associated to $\F$.
 The singular set of   $\F$ is given by the  analytic subset $$\mathrm{Sing}(\F)=\{x\in X; \omega_{\F}(x)=0 \}\subset X.$$

  For more details we refer to \cite{CMM, CCJP3}. 

\begin{defi}\index{Foliations}
A distribution  $\F$  on $X$ whose   tangent sheaf is  invariant under the Lie bracket is called a foliation, i.e, 
   if  $[T\F,T\F]\subset T\F$.
\end{defi}
Therefore,  we have  a regular foliation on $X\setminus \mathrm{Sing}(\F)$.

\subsection{Pfaff systems}
\begin{defi} \index{Pfaff systems}
A   Pfaff system $\F$  of  codimension  $p$ on $X$ is a nonzero map
$\eta : \sL^* \to  \Omega_X^p $, where  $\sL$ is a
line bundle  on $X$.
It corresponds to a  twisted $p$-form  $\omega_{\eta}\in H^0(X, \Omega_X^p\otimes\sL)$.
The  singular set  $\mathrm{Sing}(\eta)$ of $\eta$
is the zero set  of $\omega_{\eta}$. 
\end{defi}

The scheme structure of $\mathrm{Sing}(\F)$ is defined as follows. Consider the the dual morphism $\eta^* :\wedge^pTX\to \sL $. Twisting by $\sL^*$ we obtain a morphism  
$\wedge^pTX\otimes \sL^* \longrightarrow \mathcal{O}_X$. The  image   is the ideal sheaf $\iz$ of a subscheme $Z$ whose support is $ \mathrm{Sing}(\F)$.

Consider  a Pfaff system $\F$ of codimension  $p$ on $X$ induced by a twisted $p$-form $\omega \in H^0(X, \Omega_X^p\otimes\sL)$.  Then $\omega$ is determined by the following: 
\begin{itemize}
\item [(i)] an open covering $\{U_{\alpha}\}_{\alpha\in \Lambda}$  of $X$;
\item [(ii)]\   holomorphic $p$-forms $\omega_{\alpha} \in \Omega^p_{U_{\alpha}}$ satisfying
\end{itemize}
\begin{eqnarray}\nonumber
\omega_{\alpha } =  h_{\alpha \beta} \omega_{\beta}\,\,\,\,\,\,\,\,\,\,\,\,\,\,\,\,\,\,\mbox{on}\,\,\,\,\,\,\,\,\,\,\,\,\,\,\, U_{\alpha}\cap U_{\beta}\neq \emptyset,
\end{eqnarray}
\noindent where $h_{\alpha \beta}\in \OO(U_{\alpha}\cap U_{\beta})^{\ast}$ determines a cocycle representing $\sL$.

\begin{defi}\index{Invariant subvarieties}
We say that an analytic subvariety $V \subset X$ is  invariant by a Pfaff system $\F$ 
 induced by a twisted $p$-form
$\omega \in H^0(X, \Omega_X^p\otimes\sL)$ if 
$
i^*\omega \equiv 0,
$
where $i: V\hookrightarrow X$ is the inclusion map.
\end{defi}

Let   $\F$ be a Pfaff system of  codimension  $p$  on $X$ induced by a twisted $p$-form $\omega \in H^0(X, \Omega_X^p\otimes\sL)$  and $V$ an analytic subvariety of  $X$ of pure  codimension $k\leq p$ . Suppose that for each $\alpha\in \Lambda$ we have 
$$
V\cap U_{\alpha} =\{z\in U_{\alpha}: f_{\alpha,1}(z)=\cdots= f_{\alpha,k}(z) = 0\},
$$
where $f_{\alpha,1},\ldots, f_{\alpha,k}\in \OO(U_{\alpha})$. If $V$ is invariant by $\omega$, then for each $i \in \{1,\ldots,k\}$  there exist holomorphic  $(p+1)$-forms $\theta^{\alpha}_{i1},\ldots, \theta^{\alpha}_{ik}\in \Omega^{p+1}_{U_{\alpha}}$, such that
\begin{eqnarray}\label{3p2}
\omega_{\alpha}\wedge df_{\alpha,i} = f_{\alpha,1}\theta^{\alpha}_{i1} + \cdots + f_{\alpha,k}\theta^{\alpha}_{ik}.
\end{eqnarray}

A. G. Aleksandrov \index{Aleksandrov} in  \cite{Ale}  introduced  the concept of   multiple residues  of a logarithmic differential form with poles along a complete intersection which is a   generalization of   Saito's  residues \cite{Saito} \index{Saito}. Here we recall such theory.  
 Let $U$ be a germ of $n$-dimensional complex manifold and $D$  an analytic reduced hypersurface in  $U$ whose decomposition into irreducible components is  given by 
$$
D = D_1\cup\cdots \cup D_k.
$$
Suppose that the analytic subvariety $V = D_1\cap\cdots \cap D_k$  is reduced and has pure codimension $k$. We assume  that  
$$
V  = \{z\in U : f_{ 1}(z) =\cdots =f_{ k}(z) = 0\},
$$
 with $f_{1},\ldots,f_{  k}\in \OO(U )$ and for each $i \in \{1,\ldots, k \}$,
$
D_i   = \{z\in U : f_{i}(z) = 0\}.
$
Since $V$ is  a  reduced variety, the   $k$-form $df_{1}\wedge \ldots \wedge df_{ k}$
 does not vanish identically at each irreducible  component of $V$. Denote by $\Omega^q_U(\hat{D}_i)$, $q\geq 1$, the $\mathcal{O}_U$-module of meromorphic differential $q$-forms  with simple poles on  $\hat{D}_i = D_1\cup \cdots \cup  D_{i-1}\cup D_{i+1} \cup \cdots\cup D_k$, for each $i = 1,2,\ldots,k$.

We can prove the following result as a consequence of Alexandrov's theory.

\begin{prop}\label{prop_dec}   \cite[Proposition 2.8]{Correa-Diogo-GSV}
Let   $\F$ be a Pfaff system of  codimension  $p$  on a complex manifold  $X$ induced by a twisted $p$-form $\omega \in H^0(X, \Omega_X^p\otimes\sL)$,
and  consider  $V\subset X$ a reduced local complete intersection subvariety  of codimension $k$ which is  invariant by $\omega$.
Then for all local representations  $\omega_{\alpha}= \omega|_{U_\alpha}$  of $\omega$,  and all local expressions  of $V$ in $U_\alpha$
$$
V\cap U_{\alpha} = \{z\in U_{\alpha}: f_{\alpha,1}(z)=\cdots= f_{\alpha,k}(z) = 0 \},
$$
 there exist a holomorphic function $g_{\alpha}\in\OO(U_{\alpha})$, a holomorphic $(p-k)$-form $\xi_{\alpha}\in\Omega^{p-k}_{U_{\alpha}}$ and a holomorphic $p$-form $\eta_{\alpha}\in\Omega^{p}_{U_{\alpha}}$, such that
\begin{eqnarray}\label{expr1110}
g_{\alpha}\, \omega_{\alpha} = df_{\alpha,1}\wedge\cdots \wedge df_{\alpha,k} \wedge  \xi_{\alpha} + \eta_{\alpha}.
\end{eqnarray}
\noindent Moreover, $g_{\alpha}$  is not identically zero on every irreducible component of  $V$ and $\eta_{\alpha}$ is given by
\begin{eqnarray}\nonumber
\eta_{\alpha} = f_{\alpha,1}\,\eta_{\alpha,1}+\cdots+f_{\alpha,k}\,\eta_{\alpha,k},
\end{eqnarray}
\noindent where each $\eta_{\alpha,i}\in\Omega^{p}_{U_{\alpha}}$ is a holomorphic $p$-form.

\end{prop}

\begin{defi}\label{intpri}\index{First integral}
If $\F$ is a Pfaff  system  on   $X$, a first integral for $\F$ is a non-constant meromorphic function   $f : X  \dashrightarrow \mathbb{P}^1$,  such that the fibers of $f$ are invariant by $\F$.
\end{defi}

By using the isomorphism  $ \Omega_X^p \simeq \wedge^{n-p}TX \sl \otimes \det TX^*$, we have that  a  Pfaff system $\omega\in \mathrm H^0(X,\Omega_X^p\otimes \sL)$ induces a global section 
$$ \vartheta_{\omega} \in H^0(\wedge^{n-p}TX \otimes \det TX^* \otimes \sL )
$$
which is the so-called   {\it Pfaff field}  associated to $\omega$, see \cite[Section 3]{Esteves-Kleiman-Pfaff}.

\subsection{Pfaff systems on complex projective spaces}

 Let $\omega \in \mathrm H^0(\PP^n,\Omega_{\PP^n}^k(r))$ be a holomorphic  Pfaff system of codimension  $k$ on  $\PP^n$. 
Take a    generic non-invariant  linearly embedded subspace $i:H\simeq \PP^k \hookrightarrow \PP^n$.  We have an induced  non-trivial section  
$$i^*\omega \in \mathrm H^0(H ,\Omega_{H}^k(r)) \simeq \mathrm H^0( \PP^k , \OO_{\PP^k}(-k-1+r)),$$ 
since $\Omega_{\PP^k}^k=\OO_{\PP^k}(-k-1)$.  The tangency set between $\omega$ and $H$, denoted by $Z(i^*\omega)$, is defined as the hypersurface of zeros of $i^*\omega$ on $H$.  The {\it degree} \index{Degree of foliations} of $\omega$, denoted by $\deg(\omega)$, is defined as the degree of $Z(i^*\omega)$ in $H$ and, therefore, is given by $$\deg(\omega)=-k-1+r.$$
In particular,  we have that
\begin{equation}\label{secao-campo}
\omega \in \mathrm H^0(\PP^n,\Omega_{\PP^n}^k(d+k+1))\simeq H^0(\PP^n,\wedge^{n-k} T\PP^n(d-k+n)),
\end{equation}
 where $\deg(\omega) = d$. A Pfaff system of degree $d$ can be induced  by  a polynomial differential  $k$-form on $\C^{n+1}$ with homogeneous coefficients of degree $d+1$, see for instance \cite{CMM,CMM2,Jouanolou2}.

\section{Characteristic classes and residues}
\index{Characteristic classes and residues}

\subsection{Baum-Bott Theorem}

In \cite{Baum-Bott}  P. Baum and R. Bott  developed  a   residue theory
for singular holomorphic foliations on complex manifolds. More precisely, they proved the following:

\begin{teo}[Baum-Bott] \index{Baum-Bott Theorem} Let $\sF$ be a holomorphic foliation of codimension $k$ on a complex manifold $X$ and $\varphi$ be a homogeneous symmetric polynomial 
of degree $d$ satisfying $ k < d \leq n$. Let $Z$ be a compact connected component of the singular set $\mathrm{Sing}(\sF)$. Then  there exists a homology class $\mathrm{Res}_{\varphi}(\sF, Z) \in \H_{2(n - d)}(Z; \mathbb{C})$ such that:
\begin{enumerate}
\item[$i)$] $\mathrm{Res}_{\varphi}(\sF, Z) $ depends only on $\varphi$ and on the local behavior of the leaves of $\sF$ near $Z$. 
\item[$ii)$] Suppose that  $X$ is compact and denote by $\mathrm{Res}({\varphi}, \sF, Z):=\alpha_{\ast}\mathrm{Res}_{\varphi}(\sF, Z)$,
where $\alpha_{\ast}$ is the composition of the maps
$$\displaystyle \H_{2(n - d )}(Z; \mathbb{C}) \stackrel{i^{\ast}}\longrightarrow \H_{2(n - d )}(X; \mathbb{C}) $$ and $$\displaystyle \H_{2(n - d )}(X; \mathbb{C}) \stackrel{P}\longrightarrow \H^{2d}(X; \mathbb{C}) $$
 where  $i^{\ast}$ is the induced map of the inclusion $i : Z \longrightarrow X$ and $P$ is the Poincar\'e duality. Then 
 $$\displaystyle \varphi (N_{\sF}) = \sum_{Z} \mathrm{Res}({\varphi}, \sF, Z).$$
 \end{enumerate}

\end{teo}

An explicit calculation of the residues is difficult in general, see \cite{Baum-Bott,Cenkl, Bracci-Suwa, perrone,  Correa-Lourenco,Correa-Arturo, Brasselet-Correa-Lourenco,Correa-Arturo-Renato-Gilcione}.
If the foliation $\sF$ has  dimension one with isolated singularities,  Baum and Bott  in \cite{BB1}  show that residues can be expressed in terms of a  Grothendieck  residue\index{Grothendieck residue}, i.e, for each $x\in \mathrm{Sing}(\sF)$ we have 
$$
\mathrm{Res}_{\varphi}(\sF, x)=  \mbox{Res}_{x}\Big[\varphi(Jv)\frac{dz_{1}\wedge\cdots   \wedge dz_{n}}{v_{1}\cdots   v_{n}}\Big] ,
$$
where $v=\sum_{i=1}^nv_i\frac{\partial}{\partial z_i}$ is a germ  of holomorphic vector field at $x$ tangent to $\sF$ and $Jv$ is the jacobian of $v$. In particular, if  $\varphi=\det$,  then 
$$
\mathrm{Res}_{\det}(\sF, x)=  \mbox{Res}_{x}\Big[ \det(Jv)\frac{dz_{1}\wedge\cdots   \wedge dz_{n}}{v_{1}\cdots   v_{n}}\Big]= \mu_p(v) ,
$$
where $\mu_x(v)$ is the Milnor number \index{Milnor number} of $v$ at $x$. In this situation Baum-Bott theorem is 
$$
\int_Xc_n(N_{\F})=\int_Xc_n(TX-T\F)=\sum_{x\in \mathrm{Sing}(\sF) } \mu_x(v). 
$$
A similar formula is   given in \cite{Correa-Miguel-Soares} in 
the context of complex compact orbifolds and these residue formulae have been applied  to  the Poincar\'e problem for quasi-homogeneous  vector fields \cite{BCR, Correa-Soares-pesos}.

In  the  context of a  hypersurface $V$ invariant by a  one-dimensional foliation $\F$, we have the following Baum-Bott type theorem \cite{Correa-Machado-TAMS, Correa-Machado-Log}: 
if  $\F$ has  isolated singularities and $V$ is a  normal crossing divisor\index{Normal crossing divisor}, then
\begin{equation}\label{BB-log}
    \displaystyle\int_{X}c_{n}(TX(-\log\, V)- T\F) =    \sum_{x\in  \mathrm{Sing}(\F)\cap (X\setminus  V)} \mu_x(\F) + \sum_{x\in \mathrm{Sing}(\F) \cap  V }   \mathrm{ Log}(\F,V,x),
\end{equation}
where $ \mathrm{Log}(\F,V,x)$ denotes    Aleksandrov's logarithmic index \cite{Aleksandrov} of $\F$ at the point $x$ and $TX(-\log\, V)$ \index{Logarithmic tangent bundle} is the logarithmic vector bundle associated to $V$.    

\subsection{GSV-index for  holomorphic Pfaff Systems}

The GSV-index for vector fields tangent to hypersurfaces with isolated singularities  was introduced by X. G\'omez-Mont, J. Seade and A. Verjovsky \cite{GSV} 
 as a generalization of the the Poincar\'e-Hopf index.\index{Poincar\'e-Hopf index} 
The concept of GSV-index for vector fields tangent to complete intersections has been  extended   by J. Seade and T. Suwa in  \cite{SeaSuw1, SeaSuw2} and   J. -P. Brasselet, J. Seade and T. Suwa in  \cite{a03} \index{GSV-index}.

M. Brunella in \cite{Brunella}  \index{Brunella} introduced the GSV-index for one-dimensional singular foliations in complex surfaces in terms of the germs of  $1$-forms inducing the foliation and   established a relation  between  the GSV-index  with the 
  Khanedani-Suwa variational   \cite{KS} and   Camacho-Sad    \cite{CS} indices. 
Let us recall   Brunella's definition of GSV index.
 Let $X$ be a complex compact  surface and  $\F$ a  one-dimensional holomorphic foliation on $X$. Let $C$ be  a reduced curve on $X$ invariant by $\F$. 
Consider   $ \omega \in \mathrm H^0(X,\Omega_X^1\otimes   \sL)$ a twisted 1-form   inducing  the foliation  $\F$.
Given a point $x\in C$, let $f=0$  be a local equation of $C$ in a neighborhood $U_{\alpha}$ of $x$ and let $\omega_{\alpha}$ be the holomorphic $1$-form inducing the foliation   $\F$ on $U_{\alpha}$. Since $\F$ is logarithmic along $C$,  it follows from   \cite{Saito, AlcidesLog, SuwaLog} that there are holomorphic functions $g$ and $\xi$ defined in a neighborhood of $x$,  
that do not vanish identically  
simultaneously on $C$, such that
\begin{eqnarray} \label{2p11}
g\,\frac{\omega_{\alpha}}{f} = \xi\, \frac{df}{f} + \eta,
\end{eqnarray}
\noindent with $\eta$ being a suitable holomorphic $1$-form. In view of this relation, M. Brunella in  \cite{Brunella} showed that   the GSV-index can be defined as follows:
$$
\mathrm{GSV}(\F, C, x) = \sum_i \mathrm{ord}_x \left(\displaystyle\frac{\xi}{g}|_{C_i}\right),
$$
 where $C_i\subset C$ are irreducible components of $C$ and  $\mathrm{ord}_x \left(\displaystyle\frac{\xi}{g}|_{C_i}\right)$ denotes the  order of vanishing of $\displaystyle\frac{\xi}{g}|_{C_i}$ at  $x$.  
\begin{teo} [   \cite{Brunella, Brunella2}] Let $\F$ be a  one-dimensional holomorphic foliation on a complex compact surface $X$    and   $C \subset X$ a reduced curve invariant by $\F$. Then
\begin{eqnarray}\nonumber
\sum_{x\in Sing\left(\F\right) \cap C} \mathrm{GSV}(\F,C,x) =  \sL \cdot  C - C \cdot C\,.
\end{eqnarray}
\end{teo}

Cavalier  and  Lehmann in \cite{Cavalier-Lehmann} have  studied  a GSV   index via K-theory for  locally complete intersection invariant curves  and they prove certain  inequalities in the  Poincar\'e problem context. Recently, 
T. Suwa in \cite{Suwa-GSV}  gave a new interpretation of the  GSV-index as  a  residue arising from a certain localization of the Chern class of the ambient tangent
bundle.

In   \cite{Correa-Diogo-GSV}  we  introduce a  GSV  type index for Pfaff systems on projective  manifolds and we  prove  some of its important  properties. 
Consider a Pfaff system of codimension $p$ induced by the twisted $p$-form   $\omega\in \mathrm H^0(X,\Omega_X^p\otimes \sL)$  and $V$ a reduced  local complete intersection  subvariety   of pure  codimension $k$ invariant by $\omega$. Let us denote  $\mathrm{Sing}(\omega,V): =\mathrm{Sing}(\omega)\cap V  $. 
We   also assume that the codimension of the system $\omega$ coincides with the codimension of $V$, i.e., $p=k$.
Fixed an irreducible component $S_i$ of  $\mathrm{Sing}(\omega,V)$,   
let $\omega_{\alpha}= \omega|_{U_\alpha}$ be a local representative of  $\omega$,  such that $U_{\alpha}\cap S_i \neq \emptyset$.    Suppose V is represented in $U_\alpha$ by
$$
V\cap U_{\alpha} = \{z\in U_{\alpha}: f_{\alpha,1}(z)=\cdots= f_{\alpha,k}(z) = 0 \}.
$$
\noindent and  it follows from Proposition \ref{prop_dec}   that   
\begin{eqnarray}\label{2p10}
g_{\alpha}\,\omega_{\alpha} = (df_{\alpha,1}\wedge\cdots \wedge df_{\alpha,k})\, \xi_{\alpha} + \eta_{\alpha},
\end{eqnarray}
\noindent where $\eta_{\alpha} = f_{\alpha,1}\,\eta_{\alpha,1} + \cdots +f_{\alpha,k}\,\eta_{\alpha,k}$,   $\eta_{\alpha,i}\in \Omega^{k}_{U_{\alpha}}$, and   $\xi_{\alpha}$ being a holomorphic function. Now, we are  able to define 
the GSV-index for Pfaff systems with  an invariant reduced  local complete intersection  subvariety.   

\begin{defi}Suppose that $S:=\mathrm{Sing}(\omega,V)$ is a codimension one subvariety of $V$.  The GSV-index of $\omega$ relative to $V$ in $S$ is given   by 
\begin{eqnarray}\nonumber
\mathrm{GSV}(\omega, V, S) : = \sum_j \mathrm{ord}_{S}\left(\displaystyle\frac{\xi_{\alpha}}{g_{\alpha}}|_{V_j}\right),
\end{eqnarray}
where the sum is taken over all irreducible components $V_j $ of $V$ and 
 $\mathrm{ord}_{S}\left(\displaystyle\frac{\xi_{\alpha}}{g_{\alpha}}|_{V_j}\right)$ denotes  the order of vanishing of $\frac{\xi_{\alpha}}{g_{\alpha}}|_{V_j}$ along   $S$.
 \end{defi}

In  \cite{Correa-Diogo-GSV} we   prove  the following result.

\begin{teo}[ \cite{Correa-Diogo-GSV}]\label{prop9}
Let $X$ be a projective manifold and  $V\subset X$ a reduced local complete intersection subvariety  of codimension $k$ invariant by  a Pfaff system of  codimension  $k$  induced by a twisted $k$-form $\omega\in \mathrm H^0(X,\Omega_X^k\otimes \sL)$.  Then the following hold:
\begin{itemize}
\item[(a)]\  
there exists a  complex number $\mathrm{GSV}(\omega, V, S_i)$ which  depends only on the local   representatives of  $\omega, V$ and  $S_i$;
\\ 
\item[(b)] \  
if     $\mathrm{Sing}(\omega,V): =\mathrm{Sing}(\omega)\cap V$ has codimension one in $V$, then  
$$
\sum_i \mathrm{GSV}(\omega, V, S_i) [S_i] = c_1([\sL \otimes \det (N_{V/X})^{-1}])|_V\frown [V],  
$$
where $S_i$ denotes an irreducible component of $\mathrm{Sing}(\omega,V)$ and  $N_{V/X}$ is the normal sheaf of   $V$.
\end{itemize}
\end{teo}

\section{Polar varieties and  foliations}

 In   \cite{Soares-Polar}  Soares proves several bounds relating polar classes of smooth invariant varieties and the degree
of the foliation. In what follows we shall introduce the polar classes associated to a smooth projective
variety.

If  $\F$  is  a one-dimensional  foliation  on
$\mathbb{P}^{n}$   of
degree $d $, then  $\F$  is given by  a global section of $T\mathbb{P}^{n}(d-1)$, see  (\ref{secao-campo}).  
Let $D$ be an  analytic hypersurface  on  $\mathbb{P}^{n}$ defined locally by
functions $\{f_{\alpha}\in
\mathcal{O}(U_{\alpha})\}_{\alpha \in\Lambda}$, where
$\{U_{\alpha}\}_{\in\Lambda}$ is an open covering of  $\mathbb{P}^{n}$.
If
$U_{\alpha\beta}:=U_{\alpha}\cap U_{\beta}\neq\emptyset$, 
then there  exists $f_{\alpha\beta}\in
\mathcal{O}^*(U_{\alpha\beta})$, such that
$f_{\alpha}=f_{\alpha\beta}f_{\beta}$. Since  $\F$ is given by a section of  $T\mathbb{P}^{n}(d-1)$, then we have  collections
$(\{\vartheta_{\alpha}\};\{U_{\alpha}\};
\{g_{\alpha\beta}\in \mathcal{O}^*(U_{\alpha}) \})_{\alpha
\in \Lambda}$ on $\mathbb{P}^{n}$, where $g_{\alpha\beta}$ is the cocycle
inducing $ T\F^*=\mathcal{O}_{\mathbb{P}^{n}}(d-1)$. Consider the following functions
$$\zeta_{\alpha}^{(\F,D)}=\vartheta_{\alpha}(f_{\alpha})_{|D}\in \mathcal{O}(U_{\alpha}\cap D).$$
If $U_{\alpha}\cap U_{\beta}\cap D\neq\emptyset$, 
using  Leibniz's rule we get
$\zeta_{\alpha}^{(\F,D)}=f_{\alpha\beta}g_{\alpha\beta}\zeta_{\beta}^{(\F,D)}$.
This yields a global section $\zeta^{(\F,D)}$ of the line
bundle $(T\F^*  \otimes \mathcal{O}_{\mathbb{P}^{n}}(D))_{\mid D}$. The
\emph{tangency variety } of $\F$
  with $D$ is given by
  $$
\mathcal{T}(\F,D)=\{p\in D; \zeta^{(\F,D)}(p)=0\}.
  $$

\begin{defi}
Consider a pencil
of hyperplanes $\mathcal{H}=\{H_{\lambda}\}_{_{\lambda\in
\mathbb{P}^{1}}}$, with base locus  $\bigcap_{\lambda\in
\mathbb{P}^{1}}H_{\lambda}=\mathbb{L}^{n-2}$, where
$\mathbb{L}^{n-2}$ is a linear subspace of dimension  $n-2$ which
is not contained in $\mathrm{Sing}(\F)$. The \emph{ polar divisor } \index{Polar varieties}
of $\F$ with respect $\mathcal{H}$ is
$$
\mathcal{D}_{\mathcal{H}}=\bigcup_{\lambda\in \mathbb{P}^{1}
}\mathcal{T}(\F,H_{\lambda}).
$$
\end{defi}

\begin{lema} \cite[Lemma 2.1]{Soares-Polar}\label{lema-grau}
$\mathcal{D}_{\mathcal{H}}$   is a hypersurface of  degree
$ d+1$. 
\end{lema}

Let $ V\subset \mathbb{P}^{n}$ be a smooth projective variety, of dimension $n-k$, and $\mathbb{L}^{k+j-2}$   a linear subspace of dimension $k+j-2$. Then
the $j$-th polar locus \cite{Fulton} of $V$ is defined by 
$$
\mathcal{P}_{j}(V)=\{p\in V; \dim(
\mathbf{T}_p V \cap \mathbb{L}^{k+j-2}) \geq j-1  \},
$$
for $0 \leq  j \leq n -k$. The  $j$-th class  $\varrho_j(V)$ of $V$ is the degree of the cycle $\mathcal{P}_{j}(V)$
 and can be written as 
$$
\varrho_j(V)= \sum_{i=0}^{j}(-1)^i
\binom{n-k-i+1}{j-i}\int_Vc_i(TV)c_1(\mathcal{O}_V(1))^{n-k-i}.
$$

\begin{lema} \cite[Lemma 3.1]{Soares-Polar}\label{lema-polar}
Let $\mathcal{V}$ be 
  a smooth irreducible algebraic variety of dimension $n-k$, invariant by a one-dimensional foliation
$\F$  such that $\mathcal{V}$ is  not contained in 
 $\mathrm{Sing}(\sF)$. 
Consider
a pencil of hyperplanes
$\mathcal{H} =\{H_{\lambda}\}_{\lambda\in
\mathbb{P}^{1}}$.  Then
\begin{center}
$ \mathcal{P}_{n-k}(\mathcal{V})\subset
\mathcal{D}_{\mathcal{H}}$ and $\mathcal{P}_{0}(\mathcal{V})=\mathcal{V}\not\subset
\mathcal{D}_{\mathcal{H}}$.
\end{center}
\end{lema}
Mol in \cite{Mol} introduces the notion of polar sets associated to the holomorphic distributions.

\section{Poincar\'e and Painlev\'e problems for foliations and  Pfaff systems}

In \cite{Cerv1}  Cerveau and Lins Neto have initiated the study of  Poincar\'e's problem by showing 
that if  $C$ is a curve   of degree $m$   with  singularities
of normal crossing type, invariant by  a foliation  by curves  $\F$ on $\mathbb{P}^2$, of degree $d$, then  \index{Cerveau and Lins Neto bound}
$$
m\leq d+2.
$$
Moreover, they proved that  the  equality holds if $\F$  is given by a  logarithmic $1$-form and $C$ is reducible. If $C$ is irreducible this bound follows from the following equation  \cite[Proposition 885]{Cerv1} 
\begin{equation}\label{soma-indices-CL}
  2-2g(C)=\sum_{i}\imath(\F,C_i) - m(d-1),  
\end{equation}
where $\imath(\F,C_i)$  is the sum of  the multiplicity along the local branch $C_i$ of the local vector field tangent to $\F$. In fact, $ \imath(\F,C_i)=1$ for all local branches $C_i$
of $C$ through the singular points of $\F$ on $C$. If we denote by  $s$ be the number of
nodal points of $C$ and  since for each nodal point,   we have $2$ local branches of
$C$ we have that  the contributions of the nodal points in $\sum_{i}\imath(\F,C_i)$ is precisely
$2s$.  By the genus formula  we get  
$$
 2-2g(C)=-m^2+3m+2s.
$$
Defining $\ell= \sum_{i}\imath(\F,C_i)-2s\geq 0$ we have that 
$$
-m^2+3m+2s=2-2g(C)=\sum_{i}\imath(\F,C_i) - m(d-1)=2s+\ell - m(d-1) 
$$
results in 
$
 0\leq \ell=m(d+2-m),   
$ which implies that $m\leq d+2$. 

Carnicer in \cite{Carnicer} proved that if the singularities of the foliation along $C$ are of non-dicritical type,   then $m\leq d+2$ as well.

In order to investigate these results in  higher dimensions many authors have provided   generalizations of the formula (\ref{soma-indices-CL}).
Let $\F$ be a one-dimenional foliation on a projective manifold $X$   and $C\subset X$ a curve invariant by $\F$.  
Esteves and Kleiman in \cite[Proposition 5.2]{Esteves-Kleiman-caracteristica} prove the following formula
\begin{equation}\label{formula-EK} \index{Esteves and Kleiman formula}
    2p_a(C)-2- C \cdot  T\F^*=\lambda(C)-\deg(C\cap \mathrm{Sing}(\sF)), 
\end{equation}
where  $p_a(C)$ is the  arithmetic genus of $C$, $\lambda(C)$ is the colength in the dualizing sheaf of the subsheaf generated by the K\"ahler
differentials.  Esteves and Kleiman showed this result in the context of a (possibly singular) algebraic variety   over an algebraically closed field of  arbitrary characteristic. But, in this survey we will only deal with smooth ambient spaces. The formula (\ref{formula-EK}) is a generalization of  those  obtained in \cite[Prop. 2.7, p. 659.]{Alcides-Marcio} and \cite[Prop. 2.2. p. 60 ]{Campillo-Carnicer-Fuente}.

If $C \subset \mathbb{P}^n$
 is a smooth curve which is a  complete
intersection of hypersurfaces of degrees $d_1, \dots,  d_{n-1}$  invariant by a one-dimensional  foliation $\F$ of degree $d$,   Soares has provided in  \cite{Soares-Annals}  the following bound 
\begin{equation} \index{Soares bound}
  d_1+ \cdots + d_{n-1} \leq d+n-1. 
\end{equation}
Now, if  $C$ has at most ordinary nodes  as  singularities, Campillo, Carnicer and J.  
de la Fuente showed in \cite{Campillo-Carnicer-Fuente},  by using a   formula of the type of  (\ref{formula-EK}),  that
\begin{equation}
  d_1+ \cdots + d_{n-1} \leq d+n. 
\end{equation}
Esteves in \cite{Esteves-Z} also  proved  the same bound by employing  algebro-geometric techniques.

Let $V$ be a  \index{Complete intersection} variety of codimension $k$  
which is a  complete intersection of hypersurfaces of degrees
$d_1,\dots,d_{k}$.  If $V$ is invariant  by  a Pfaff system  $\F$ of degree $d$ and   codimension 
$k$ on $\mathbb{P}^n$, Esteves and Cruz proved  in  \cite[Corollary 4.5]{Cruz-Esteves}  that 
$$
d_1+\cdots+d_{k}\leq \left\{
  \begin{array}{ll}
    d+k, & \hbox{if } \ \rho\leq 0  \\
\\
    d+k+\rho, & \hbox{if } \rho>0
  \end{array}
\right.
$$
where $\rho:= \sigma+k+1-d_1-\cdots-d_{k}$, with $\sigma$
denoting the Castelnuovo-Mumford regularity\index{Castelnuovo-Mumford regularity} of the singular locus of
$V$. In \cite{Correa-Jardim-genera} we prove that 
if $V$ is nonsingular in codimension $1$, then one can take $\rho=1$, regardless
of $\sigma$. That is, under this condition we get that
\begin{equation}\label{cota-IC}
d_1+\cdots+d_{k}\leq d+k+1.
\end{equation}
We will see in the section \ref{section-GSV} that  the non-negativity of the GSV-indices  for Pfaff systems  yields  an  obstruction to obtain  the bound (\ref{cota-IC}).

\subsection{The approach of Soares }
Soares  made substantial contributions to the investigation of the Poincar\'e problem for one-dimensional   foliations   with smooth invariant varieties. He has  obtained certain bounds via Baum-Bott Theorem and also polar varieties.  

Let $\sF$  be  a one-dimensional foliation on $\mathbb{P}^n$,  of degree $d$,  with non-degenerate  
 isolated singularities.  
If $V\subset \mathbb{P}^n$ is a smooth hypersurface of degree $m$ which is invariant by  $\sF$, then Soares  applied the Baum-Bott Theorem 
 to prove that ( \cite[Theorem A]{Soares-Inventiones} )
$$
\sum_{i=0}^{n-1}[1+(-1)^i(m-1)^{i+1}]d^{n-1-i} =\int_V c_{n-1}(TV(d-1))=\sum_{p\in \mathrm{Sing}(\sF) \cap V }  \mu_p(v)=\# \mathrm{Sing}(\sF)  \cap V.
$$
In order to prove the bound
\begin{equation}\label{Soares-cota-1}
m\leq d+1
\end{equation}
Soares   studied  the following  generating function 
$$
\psi(x)=\sum_{i=0}^{n-1}[1+(-1)^ix^{i+1}]d^{n-1-i}.
$$
By  Lehmann's result  \cite{Lehmann} we have that $\psi(m-1)=\#   \mathrm{Sing}(\sF)  \cap V$  is positive.  If  $n$ is even, then $\psi(x)\leq 0$   whenever $x\geq d+1$. Thus $m\leq d+1 $. For the  case when  $n$ is odd, then 
$$
\psi(x) \geq  \psi(m-1)=\# \mathrm{Sing}(\sF)  \cap V 
$$
whenever   $x\geq d $ and once  again we can conclude  that  $m\leq d+1 $. 

As we have seen, Soares'  approach consists of studying the number of singularities contained in the invariant hypersurface. In \cite[Theorem 3]{Correa-Machado-TAMS} the authors  by using the formula (\ref{BB-log}) prove that if  $n$ is odd, then $\mathrm{Sing}(\sF) \subset V$ if and only if the Soares' bound  (\ref{Soares-cota-1})  is achieved.

Now,      if  $V\subset \mathbb{P}^n$  
 is a smooth irreducible  projective variety which is a  complete
intersection of hypersurfaces of degrees $d_1, \dots,  d_{r}$  invariant by a foliation $\F$, of degree $d$ and with  non-degenerate isolated singularities.    Soares, by  applying Baum-Bott theorem \cite{Soares-Annals},  provided the following bound 
\begin{equation}\label{Soaresb1}
 \frac{\varrho_{n-r}(V)}{\varrho_{n-r-1}(V)} \leq d.  
\end{equation}
On the other hand, Soares proved the following bound via polar varieties    
\begin{equation}
 \frac{\varrho_{n-r-j}(V)}{\varrho_{n-r-j-1}(V)} \leq d+1,  
\end{equation}
if $ \mathcal{P}_{n-k-j}(\mathcal{V})\subset
\mathcal{D}_{\mathcal{H}}$ but $ \mathcal{P}_{n-k-j-1}(\mathcal{V})\not\subset
\mathcal{D}_{\mathcal{H}}$, for some $0\leq j\leq n-k-1$.  Observe  from Lemma \ref{lema-polar} that 
$ \mathcal{P}_{n-k}(\mathcal{V})\subset
\mathcal{D}_{\mathcal{H}}$  and  $  \mathcal{V}\not\subset
\mathcal{D}_{\mathcal{H}}$, then 
\begin{equation}\label{Soaresb2}
 \frac{\varrho_{n-r}(V)}{\varrho_{n-r-1}(V)} \leq d+1.  
\end{equation}
Therefore, under the assumption that the foliation has non-degenerate singularities the bound (\ref{Soaresb1}) is
sharper than the one in (\ref{Soaresb2}).

\subsection{Brunella-Mendes  approach  }

 In \cite{BruMe} Brunella and Mendes
  proved a bound for  normal crossing hypersurfaces invariant by Pfaff systems as a consequence of results  due to Deligne \index{Deligne} \cite{Deligne} and Saito \cite{Saito}    on  the theory of logarithmic forms. 
Brunella and Mendes show that a Pfaff system, of codimension $k$,  with an invariant hypersurface is a meromorphic  logarithmic form, i.e, a global section of $$\Omega_{\mathbb{P}^n}^1(\log V)\otimes \mathcal{O}_{\mathbb{P}^n}(-V) \otimes \sL.  $$ 

    They apply    Bogomolov's lemma \cite{Bogomolov} \index{Bogomolov} which states that if there exists a nontrivial
 global section  of $\Omega_{ X}^1(\log V)\otimes   N $, then 
\begin{equation}\label{Bogo-ine}
    \kappa(X, N^{-1} )\leq k, 
\end{equation}
 where  $  \kappa(X, N^{-1})$ denotes the  Kodaira-Iitaka dimension \index{Kodaira-Itaka dimension} of the  line bundle $N^{-1}$. Thus, if $X= \mathbb{P}^n$ and   $V$ has degree $m$,  we have that $N^{-1}= \mathcal{O}_{\mathbb{P}^n}(m) \otimes \mathcal{O}_{\mathbb{P}^n}(-d-k-1)=\mathcal{O}_{\mathbb{P}^n}(\ell)$  for some $\ell\in \mathbb{Z}$.  
  Since $\ell>0$ implies that 
  $ \kappa(\mathbb{P}^n,  \mathcal{O}_{\mathbb{P}^n}(\ell) )=n$
   we conclude from (\ref{Bogo-ine}) that 
   $$
   m\leq  d+k+1.
   $$
 
 \subsection{Poincar\'e problem and birational geometry of foliations} \index{Birational geometry of foliations}
 The birational classification  of foliations on surfaces has  been established by Brunella \cite{Brunella-MMP, Brunella-livro},  McQuillan\index{McQuillan} \cite{McQuillan}, Mendes \cite{Mendes}.

 In \cite{Pereira} Pereira has considered the  Poincar\'e problem by obtaining  bounds for the degree of the first
integral of a foliation of general type (Kodaira dimension equal to 2)  in terms of the degree, the birational invariants   and
the geometric genus of a generic leaf of the foliation. A refinement of this result was given by Pereira and Svaldi in \cite{Pereira-Svaldi}. They prove that if  $\F$  is a foliation on $\mathbb P^2$ such that  $\F$ is birationally
	equivalent to a non-isotrivial fibration of genus $g\ge 2$.
	Then  the degree of the general leaf
	of $\F$ is bounded by
$$
\Big(4 \Big(42(2g-2)\Big)! \Big)^2 (4g -4) \deg(\F).
$$

In \cite{Hacon-Langer}  Hacon and Langer prove a result on the effective generation of pluri-canonical
linear systems of general type foliations  on  surfaces which is related  to the  Poincar\'e problem. They prove  the following:
for any  integer valued function  $P:\mathbb Z_{\geq 0}\to \mathbb Z$ and any integer $g\geq 0$, there exists an integer $\delta>0$ such that if $(X,\F )$ is  a weak nef model of a complete foliated surface with     $$\chi (X,\mathcal{O}_X(mT\F^* ))=P(m)$$ for all $m\geq 0$, and  if $(X,\F)$ has a meromorphic first integral whose general fiber has geometric genus $g$, then the general leaf $C$  satisfies 
$$C\cdot T\F^*\leq \delta.$$
 Such result is also related with the problem of  bounding  the order of   the automorphism group of foliations of general type, see \cite{Correa-Muniz, Correa-Fassarella, Muniz-Rudy,Spicer-Svaldi}.
 In the recent paper  \cite{HLT} the authors prove a bound for other birational invariant of a foliation, called  the slope, which is related to the Poincar\'e problem.

\subsection{Zariski-Esteves aproach }

Esteves in \cite{Esteves-Z} \index{Zariski-Esteves}   gave an algebraic proof of Soares' result. He proved a  characterization of all vector fields that leave invariant a  smooth
hypersurface via the Koszul complex associated to the  ideal generated by the partial derivatives
of a  polynomial defining  the hypersurface. Esteves  point out in his paper  that such proof is inspired by   an observation by O. Zariski \cite[Part c of Example 7, p. 892]{Lipman-Z}.

Let  $V=\{f=0\}$ be a smooth  hypersurface    $f$ is a polynomial  of degree $m$ defining the   hypersurface $V$ and $v$ a polynomial vector field inducing a foliation on $\mathbb{P}^n$ which leaves $V$ invariant.
Esteves in \cite{Esteves-Z} proves that 
\begin{equation}
v=\sum_{i<j}P_{i,j}\left(\frac{\d f}{\d z_i} \frac{\d }{\d z_j} -\frac{\d f}{\d z_j} \frac{\d }{\d z_i}\right) +\frac{h}{k}\sum_{i=0}^nz_i\frac{\d }{\d z_i} ,
\end{equation}
with   $P_{i,j},h\in \mathbb{C}[z_0,\dots,z_n]$. 
Now, Soares' bound is a consequence of this normal form. Indeed,  since $v\neq 0$, we conclude that  $P_{i,j}\neq 0$ for
some $i, j$. 
 Otherwise, we would have that 
$$v= \frac{h}{k}\sum_{i=0}^nz_i\frac{\d }{\d z_i},$$
and such vector field does not define a foliation on $\mathbb{P}^n$.
Thus, since $v$ has degree $d$, we have  that  $0\leq \deg P_{i,j}=d-m+1$ which gives us that $m\leq d+1$. It is worth mentioning that Esteves  considered this problem  in arbitrary characteristic.

  \begin{exem} {\rm Let $\F$ be the foliation on $\mathbb{P}^3$ induced by the polynomial vector field
\begin{eqnarray}\nonumber
\displaystyle v &=& (-z_1^{m-1} - z_2^{m-1}- z_3^{m-1})\frac{\partial}{\partial z_0} + (z_0^{m-1} - z_2^{m-1}- z_3^{m-1})\frac{\partial}{\partial z_1} +\\\nonumber\\\nonumber &&+ (z_0^{m-1} + z_1^{m-1}- z_3^{m-1})\frac{\partial}{\partial z_2} + (z_0^{m-1} + z_1^{m-1}+ z_2^{m-1})\frac{\partial}{\partial z_3}.
\end{eqnarray}
\noindent The hypersurface $V =\{z_0^{m} +z_1^{m} +z_2^{m} +z_3^{m}=0\}$  is    invariant by $\F$ and   $\mathrm{Sing}(\sF) \subset V$.
Note that  $\deg(V) = m$ and $\deg(\F) = m-1$. }
\end{exem}

\subsection{ Poincar\'e problem and  GSV-indices }\label{section-GSV}
M. Brunella in \cite{Brunella}   showed that the non-negativity of the GSV-index  is an obstruction  to the solution of the Poincar\'e problem in complex compact surfaces. 
We also have that the non-negativity of the GSV-index gives an  obstruction to the solution of the Poincar\'e Problem for Pfaff systems on complex projective  space. More precisely, 
let $\omega \in \mathrm H^0(\PP^n,\Omega_{\PP^n}^k(d+k+1))$ be a holomorphic  Pfaff system of  codimension $k$ and degree $ d$. Let $V\subset \PP^n$  be a reduced   complete intersection subvariety, of codimension $k$ and multidegree $(d_1,\dots,d_{k})$ ,  invariant  by $\omega$. Suppose that  $\mathrm{Sing}(\omega,V)$ has codimension one in $V$,  then by Theorem \ref{prop9} we obtain 
\begin{eqnarray} \nonumber 
\sum_i \mathrm{GSV}(\omega,V,S_i)\deg(S_i) \,\,\,=\,\,\, [d+k+1 - (d_1+\cdots+d_k)]\,\cdot (d_1\cdots d_k),
\end{eqnarray}
where $S_i$ denotes an irreducible component of $\mathrm{Sing}(\omega,V)$.
Therefore,   if $\mathrm{GSV}(\omega,V,S_i)\geq 0$, for all $i$, we have 
$$
d_1+ \cdots+ d_k \leq  d+k+1.
$$
In \cite{Correa-Diogo-GSV}  we  prove the following   non-negativity properties    for    the index:
\begin{itemize}
\item[(i)]  If $S_i\cap \mathrm{Sing}(V)= \emptyset$, then $\mathrm{GSV}(\omega,V,S_i)\geq 0$. 

\item[(ii)] \ If $V$ is  smooth, then $\mathrm{GSV}(\omega,V,S_i)>0$. 
\end{itemize}

Now, we  give an optimal example.  
\begin{exem}{\rm
Consider the  Pfaff system 
  $\omega \in \mathrm H^0(\PP^n, \Omega_{\PP^n}^k(d+k+1))$
given by 
$$
\omega = \sum_{0 \le j
\le k} (-1)^j d_j f_j \ df_0 \wedge \dots \wedge \widehat
{df_j} \wedge \dots \wedge df_k,
$$
where $f_j  $ is a homogeneous polynomial
of degree $d_j$. We can see    that 
$$
d_0+d_1 +\dots + d_k = d + k + 1.
$$
Suppose that  $\deg(f_0)=d_0=1$ and that $V=\{ f_1=\cdots=f_{k }=0\}$  is smooth. We have that $V$ is  invariant by $\omega$ and  
$$
d_1 +\dots + d_k = d + k .
$$}
\end{exem}

\subsection{Lins Neto's examples}
In general the Poincar\'e  problem has no solution as we can see in the following simple example. 
 Consider  the   vector field  
$$
v_{(p,q)} =px \frac{\d }{\d x}+qy \frac{\d }{\d y} 
$$
with  $p\neq q$ integers. We have that the curve  $C_{p,q} =\{y^{p} -x^{q} =0\}$ is invariant by   $v_{(p,q)}$  whereas the  family of vector fields $\{ v_{(p,q)}\}$  has  degree one.
Such family of examples has a singularity which is of dicritical type. Under certain conditions we can give  a positive answer which can be found  for instance in the works \cite{Galindo-Monserrat-dicritico,Genzmer-Mol, Campillo-Carnicer}.
Lins Neto in \cite{Alcides-exemplo} \index{Lins Neto examples} has   provided  non-trivial examples of one-parameter families which show that Poincar\'e and Painlev\'e problems have a negative answer.
Such family  of foliations on    $\mathbb{P}^2$  has small degree  and has   algebraic leaves of arbitrarily large degree. Lins Neto's families of foliations are given, in   affine chart, by the following   polynomial vector fields
$$
v_{\alpha}^2= (4x-9x^2+y^2+\alpha(2y-4xy))\frac{\partial}{\partial x}+ (6y-12xy+3\alpha(x^2-y^2))\frac{\partial}{\partial y};
$$
$$
v_{\alpha}^3= (-x+2y^2-4x^2y+x^4-\alpha(2y-x^2+xy^2))\frac{\partial}{\partial x}+ (y(-2-3xy+x^3)-\alpha(3xy-x^3+2y^3))\frac{\partial}{\partial y};
$$
$$
v_{\alpha}^4= ((x^3-1)(x-\alpha y^2))\frac{\partial}{\partial x}+ ((y^3-1)(x-\alpha x^2))\frac{\partial}{\partial y}. 
$$
Lins Neto showed that  there exists a countable and dense set of parameters $E\subset \mathbb{P}^1$
such that for any $\alpha \in E$
the induced foliation  $\F_{\alpha}$ has a rational first integrals whose     degrees  can be chosen arbitrarily large.  
In these examples the genus of an  algebraic invariant  curve   is constant and is equal to 1.  
For other interesting examples we refer  to the papers by Picard  \cite[pg 298-299]{Picard} and Movasati \cite{Movasati}.

\section{Darboux-Jouanolou-Ghys theorem for Pfaff systems} \index{Meromorphic first integral} 
In this section we will give an  idea of the proof of the  Darboux-Jouanolou-Ghys theorem for Pfaff systems. Let $X$  be a complex manifold. We denote by $\Omega_{cl}^{1}$ the sheaf of germs of closed holomorphic 1-forms and $\mathscr{M}(X)$ is the field of meromorphic functions on $X$. 
 
\begin{teo} \cite{CMM,CMM2, Jouanolou1, Ghys,Correa-dim1} \label{DJG-Pfaff}
Let $\F$ be a Pfaff system on a compact, connected, complex
manifold $X$,  induced by an $r$-form $\omega\in H^{0}(X,
\Omega^{r} \otimes \mathscr{L})$.  If $\F$ admits
\begin{equation}\label{G0}
\dim_{\mathbb{C}}H^{1}(X, \Omega_{cl}^{1}) + \dim_{\mathbb{C}}\left(H^{0}(X,
\Omega^{r+1} \otimes \mathscr{L}) / \omega \wedge H^{0}(X,
\Omega_{cl}^1)\right) + r + 1
\end{equation}
invariant irreducible analytic hypersurfaces, then $\F$ admits a
meromorphic first integral.
 \end{teo}
 
 Denote by $\mathrm{Div}(X,\F)$ the abelian group of divisors on $X$ which are
invariant by $\F$. We have the homomorphism
\begin{equation}\label{G1}
\begin{array}{c}
 \mathrm{Div}(X,\F) \longrightarrow \mathrm{Pic}(X) \\
  \\
\hskip 60pt \sum\limits_\alpha\, \lambda^\alpha\, L^\alpha \longmapsto \bigotimes\limits_\alpha\, [L^\alpha]^{\otimes \lambda^\alpha}, \lambda^\alpha \in \mathbb{Z}.
\end{array}
\end{equation}
Since $\mathrm{Pic}(X) \simeq H^1(X, \mathcal{O}_X^\ast)$, logarithmic differentiation defines a homomorphism
\begin{equation}\label{G2}
\begin{array}{c}
H^1(X,  \mathcal{O}_X^\ast) \longrightarrow H^{1}(X, \Omega_{cl}^1) \\
  \\
\hskip 34pt g \longmapsto \displaystyle\frac{d g}{g}.
\end{array}
\end{equation}
Composition of (\ref{G1}) and (\ref{G2}) gives a
$\mathbb{C}$-linear map
\begin{equation}\label{G3}
\Psi: \mathrm{Div}(X,\F) \otimes \mathbb{C} \longrightarrow H^{1}(X, \Omega_{cl}^1)
\end{equation}
which is expressed, in terms of a sufficiently fine open cover $\{U_i\}_{i \in \Lambda}$ of $X$ by: if $L^\alpha$ is defined by $f_i^\alpha =0$ in $U_i$ and, in $U_{i j}= U_i \cap U_j$, $f_i^\alpha = g_{i j}^\alpha \, f_j^\alpha$,
\begin{equation}\label{G4}
\Psi\left(\sum\limits_\alpha\, \lambda^\alpha \, L^\alpha\right)  = \left[\sum\limits_\alpha\, \lambda^\alpha \, \frac{dg_{i j}^\alpha}{g_{i j}^\alpha}\right].
\end{equation}

Consider the kernel of $\Psi$. Thus,  $\sum\limits_\alpha\, \lambda^\alpha \, L^\alpha \in \ker \Psi$ implies that  there are closed holomorphic 1-forms $\varpi_i$ such that in $U_{i j}$,
\begin{equation}\label{G5}
\sum\limits_\alpha\, \lambda^\alpha \, \frac{dg_{i j}^\alpha}{g_{i j}^\alpha}=
\varpi_j - \varpi_i.
\end{equation}
But this says that
\begin{equation}\label{G6}
\sum\limits_\alpha\, \lambda^\alpha \, \frac{df_{i}^\alpha}{f_{i}^\alpha} + \varpi_i = \sum\limits_\alpha\, \lambda^\alpha \, \frac{df_{j}^\alpha}{f_{j}^\alpha} + \varpi_j.
\end{equation}
These glue together to give a global closed meromorphic 1-form $\eta$ on $X$, defined up to  addition of a global closed holomorphic 1-form $\rho$.

Since $L^\alpha$ is $\omega$-invariant and hence $(\omega \wedge  df_{i}^\alpha)_{\mid(f_{i}^\alpha=0)} \equiv 0$,  $\omega \wedge \eta$ is a holomorphic $r+1$-form, defined up to addition  of $\omega \wedge \rho$, with $\rho$ a global closed holomorphic 1-form. Therefore, the   $\mathbb{C}$-linear map
\begin{equation}\label{G7}
\begin{array}{ccc}
\Theta : \ker (\Psi) & \longrightarrow & \mathrm{H}^{0}(X, \Omega_X^{r+1} \otimes \mathscr{L})/ \omega \wedge \mathrm{H}^{0}(X, \Omega_{cl}^1) \\
 \\
\sum\limits_\alpha\, \lambda^\alpha \,L^\alpha & \longmapsto & \overline{\omega \wedge \left(\eta + \rho\right)} \hskip 100pt
\end{array}
\end{equation}
is well-defined.

\begin{lema}\cite[Lemme 3.1.1, p. 102]{Jouanolou2}\label{Jou}
Let  $\mathscr{M}^1$  be the sheaf of germs of meromorphic 1-forms on $X$ and   $\mathscr{Q}_{X}^1:=\mathscr{M}^1/\Omega_X^1$. Then, 
the $\mathbb{C}$-linear map
\begin{equation}\label{G9}
\begin{array}{ccc}
\mathrm{Div}(X, \F)\otimes \mathbb{C} & \longrightarrow & \mathscr{Q}_{X}^1  \\
  \\
\displaystyle\sum_{\alpha}\lambda^\alpha \cdot L^\alpha & \longmapsto & \displaystyle\sum_{\alpha}\lambda^\alpha\frac{df^{\alpha}_i}{f^{\alpha}_i}
\end{array}
\end{equation}
is injective provided the divisors have no common factor.
\end{lema}

Suppose  that   $\F$ admits at least
\begin{equation}\label{G12}
\dim_{\mathbb{C}}H^{1}(X, \Omega_{cl}^1) + \dim_{\mathbb{C}}\left(H^{0}(X,
\Omega_X^{r+1} \otimes \mathscr{L}) / \omega \wedge H^{0}(X,
\Omega_{cl}^1)\right) + r + 1
\end{equation}
invariant irreducible analytic hypersurfaces.
From (\ref{G3})  we have
\begin{equation}\label{G13}
\dim_{\mathbb{C}}\ker (\Psi) \geq \dim_{\mathbb{C}}(\mathrm{H}^{0}(X,
\Omega_X^{r+1} \otimes \mathscr{L}) / \omega \wedge H^{0}(X,
\Omega_{cl}^1)) + r + 1.
\end{equation}
Hence
$\dim_{\mathbb{C}} \ker (\Theta) \geq r+1$ and $\ker (\Theta)$ is non trivial. Now, given a non zero element
 $x \in \ker( \Theta)$,
we can choose $L^{1}, \ldots ,L^{k}\in \mathrm{Div}(X, \F)$ and  $\lambda^{1},
\ldots , \lambda^{k} \in \mathbb{C}^{\ast}$ such that
\begin{itemize}
\item [$i)$]$x=
\displaystyle\sum_{\alpha = 1}^{k} \lambda^{\alpha} \, L^{\alpha} \in
\ker( \Psi)\setminus \{0\} $.
  \\
\item [$ii)$]$\displaystyle \Theta \left( x\right) = \overline{\,0} \in  H^{0}(X, \Omega_X^{r+1} \otimes
\mathscr{L}) / \omega \wedge H^{0}(X, \Omega_{cl}^1).$
\end{itemize}

Using (\ref{G6}) we have that there exists  $\mu =
(\mu_i) \in \mathcal{Z}^{0}(X, \Omega_{cl}^1)$ such that, in $U_i$,
\begin{equation}\label{G14}
\omega \wedge \left(\varpi_i + \sum_{\alpha = 1}^{k}
\lambda_{\alpha}\frac{df_{i}^{\alpha}}{f_{i}^{\alpha}}\right) = \omega
\wedge \mu_{i},
\end{equation}
which amounts to
\begin{equation}\label{G15}
\omega \wedge \left(\varpi_{i} - \mu_{i}+ \sum_{\alpha =
1}^{k} \lambda_{\alpha}\frac{df_{i}^{\alpha}}{f_{i}^{\alpha}}\right) = 0
\end{equation}
in each $U_i$. Thus we get a global closed meromorphic
$1$-form $\widetilde{\xi}$ with
\begin{equation}\label{G16}
\widetilde{\xi}_{|U_i}=\varpi_{i} - \mu_{i}+ \sum_{\alpha =
1}^{k} \lambda_{\alpha}\frac{df_{i}^{\alpha}}{f_{i}^{\alpha}}
\end{equation}
such that
\begin{equation}\label{G17}
\omega \wedge
\widetilde{\xi} = 0.
\end{equation}
In this way we can construct   $r+1>1$ global closed meromorphic $1$-forms $\widetilde{\xi}_{1}, \ldots
,\widetilde{\xi}_{r+1}$ on  $X$ such that
\begin{equation}\label{G18}
\omega \wedge
\widetilde{\xi}_{j} = 0\,, \qquad 1 \leq j \leq r+1.
\end{equation}
Define $\alpha_{1}=\widetilde{\xi}_{1} \wedge \cdots \wedge
\widetilde{\xi}_{r}$ and $\alpha_{2}=\widetilde{\xi}_{2} \wedge \cdots \wedge \widetilde{\xi}_{r+1}$ and remark that $|\alpha_{1}|_{\infty}\neq |\alpha_{2}|_{\infty}$, where $|\cdot|_{\infty}$ denotes the set of poles.
We have the following possibilities:
\\
\begin{enumerate}
  \item []\textbf{ Case 1}. $\alpha_{1} \neq 0$ and
$\alpha_{2} \neq 0.$\\
  \item [] \ \textbf{Case 2}. $\alpha_{1}=0$ or
$\alpha_{2}=0.$
\\
\end{enumerate}

\noindent \textbf{Case 1:} Since
$|\alpha_{1}|_{\infty}\neq |\alpha_{2}|_{\infty}$
we have  that
$\alpha_{1}$ and $\alpha_{2}$  are linearly independent
 over  $\mathbb{C}$. We can prove that there exist $R_1,\,R_2\in \mathscr{M}(X)$
such that  $\omega = R_{i}\,\alpha_{i}$, $i=1,2$. This gives us that 
$$d\left(\displaystyle\frac{R_{1}}{R_{2}}\right) \wedge \omega=0.$$ Moreover, since
$\alpha_{1},\alpha_{2}$  are linearly independent  over  $\mathbb{C}$, then  $ \displaystyle\frac{R_{1}}{R_{2}}$ is not constant.
This shows that the meromorphic function  $ \displaystyle\frac{R_{1}}{R_{2}}$ is a first integral for $\omega$.\\

\noindent \textbf{Case 2:}  Suppose without loss of generality  that  $\alpha_{1}=0$. Let  $m$ be the
largest integer such that $\widetilde{\xi}_{1}, \ldots ,\widetilde{\xi}_{m}$ are
linearly independent over  $\mathscr{M}(X)$. Then,
\begin{equation}\label{G20}
\widetilde{\xi}_{m+1}=\sum_{i=1}^{m}R_{i} \;
\widetilde{\xi}_{i}
\end{equation}
with $R_{1},\dots, R_m \in \mathscr{M}(X)$. Since $\widetilde{\xi}_{i}$ is  closed, for  $i=1,\dots,m+1$, we
get
\begin{equation}\label{G21}
0=\sum_{i=1}^{m}dR_{i} \wedge \widetilde{\xi}_{i}.
\end{equation}
Then, for each $j=1,\dots,m$, multiplying (\ref{G21}) by $ \widetilde{\xi}_{1}
\wedge \cdots \wedge \widehat{\widetilde{\xi}_{j} }\wedge \cdots \wedge
\widetilde{\xi}_{m}$ we obtain

\begin{eqnarray*}
0&=&\sum_{i=1}^{m}dR_{i} \wedge \widetilde{\xi}_{i} \wedge \widetilde{\xi}_{1} \wedge
\cdots \wedge \widehat{\widetilde{\xi}_{j} }\wedge \cdots \wedge
\widetilde{\xi}_{m}\\
&=&(-1)^{j+1}dR_{j} \wedge \widetilde{\xi}_{1} \wedge \cdots \wedge \widetilde{\xi}_{m}.
\end{eqnarray*}

Since  $\widetilde{\xi}_{1}, \cdots ,\widetilde{\xi}_{m}$  are linearly independent
over  $\mathscr{M}(X)$, there exist $g_{1},
\cdots ,g_{m}$ $\in  \mathscr{M}(M)$ such that  $dR_{j}=\sum_{i=1}^{m}g_{i}
\; \widetilde{\xi}_{i}.$ By $(\ref{G18})$ we get
$$
dR_{j} \wedge \omega=\sum_{l=i}^{m}g_{i}
\; \widetilde{\xi}_{i} \wedge \omega=0.
$$
Moreover, from $(\ref{G20})$ and Lemma \ref{Jou}
there exists $i_{0}\in \{1,\dots,m\}$ such that
$R_{i_{0}}$ is not constant. That is, $R_{i_{0}}$ is a meromorphic
first integral for  $\omega$. This proves the theorem.

\subsection{Comments}
We observe that if  $\F$ is a codimension one distribution, defined by  $\omega\in H^{0}(X,
\Omega^{1}_X \otimes \mathscr{L})$  which   admits
$$
\dim_{\mathbb{C}}H^{1}(X, \Omega_{cl}^{1}) + \dim_{\mathbb{C}}\left(H^{0}(X,
\Omega_X^{2} \otimes \mathscr{L}) / \omega \wedge H^{0}(X,
\Omega_{cl}^1)\right) + 2
$$
invariant irreducible analytic hypersurfaces, then $\F$ has  a
meromorphic first integral, in particular $\F$ is a foliation with  infinitely many compact leaves. In general, if $\F$ has codimension $r$, then   it follows from the proof of  Theorem \ref{DJG-Pfaff} that if either  $\alpha_{1}=\widetilde{\xi}_{1} \wedge \cdots \wedge
\widetilde{\xi}_{r}\neq 0$  or  $\alpha_{2}=\widetilde{\xi}_{2} \wedge \cdots \wedge \widetilde{\xi}_{r+1}\neq 0 $, then the Pfaff system $\F$ is given by the meromorphic  
decomposable closed form  $\omega/R_{i} = \,\alpha_{i}$, where $R_{i}$ is a meromorphic function. This shows us that  the Pfaff  system $\F$ is integrable and it is given by a global decomposable  meromorphic differential $r$-form  which is a wedge product of closed logarithmic $1$-forms. 

In \cite{Go} Gomez--Mont provides  an extension
of a theorem due to R. D. Edwards, K. C. Millett and D. Sullivan
concerning foliations with all leaves compact. He proved that if a
singular holomorphic foliation $\F$, of codimension $q$, on a
projective manifold  have  all leaves algebraic, then $\F$ is
tangent to the fibers  of  a  rational map.

In \cite{Pereira-Spicer} Pereira and Spicer   prove  an  improvement for the Darboux-Jouanolou theorem in the codimension  one case. They proved that if  $\F$ is  a codimension one foliation on a projective  manifold $X$,  $s$ is the number of   hypersurfaces invariant by $\F$ and $$\ell:=\dim NS(X) + \dim_{\mathbb{C}}  H^0(X,T\F^*) - \dim_{\mathbb{C}}  H^0(X,\Omega^1_X),$$
then the following assertions hold true.
\begin{enumerate}
\item[i)] If  $s \ge \ell$ then $\F$ is transversely affine.
\item[ii)] If $s \ge \ell +1$ then $\F$ is defined by a closed logarithmic $1$-form.
\item[iii)] If $s \ge  \ell+2$ then $\F$  has a rational first integral.   
\end{enumerate}
Here $NS(X)$ denotes the  Neron-Severi group of $X$.
 In the same work, Pereira and Spicer
show  a    Darboux-Jouanolou-Ghys type  Theorem by  providing  a characterization of codimension  one  foliations which are pull-backs
of foliations on surfaces by rational maps. 

We refer to  \cite{Scardua}
 for  general results on transversely affine and projective structures for foliations\index{Projective structure for foliation}, where Sc\'ardua also has considered   the Poincar\'e problem for codimension one foliations, see  \cite[Theorem 4.2 and 4.3]{Scardua}.

\bibliographystyle{abbrvurl}


\end{document}